\documentclass[11pt,psfig]{article}
\usepackage{amsmath,amssymb,epsfig}

\newtheorem{theorem}{Theorem}[section]

\newtheorem{lemma}[theorem]{Lemma}

\newtheorem{remark}[theorem]{Remark}

\def\sF{{\cal F}}

\def\sF{{\cal F}}

\makeatletter 
\@addtoreset{equation}{section}

\makeatother 
\begin{document}
\begin{center}
{\Huge\bf  A Portfolio Decomposition Formula}
\mbox{}\\
\vspace{2cm}
Traian A. Pirvu\\
Department of Mathematics\\
The University of British Columbia\\
Vancouver, BC, V6T1Z2\\
tpirvu@math.ubc.ca\\
\vspace{1cm}
Ulrich G. Haussmann\footnote{Work
supported by NSERC under research grant 88051 and NCE grant 30354 (MITACS).}\\
Department of Mathematics\\
The University of British Columbia\\
Vancouver, BC, V6T1Z2\\
uhaus@math.ubc.ca\\

\mbox{}\\

\today
\end{center}

\noindent {\bf Abstract.}
This paper derives a portfolio decomposition formula when the agent
maximizes utility of her wealth at some finite planning horizon. The financial market is complete and consists of multiple risky assets (stocks)
plus a risk free asset. The stocks are modelled as exponential Brownian motions with drift and volatility being It\^{o} processes. The optimal portfolio has
two components: a myopic component and a hedging one. We show that the myopic component is robust with respect to stopping
times. We employ the Clark-Haussmann formula to derive portfolio's hedging component.

\vspace{1cm}


\begin{quote}

\end{quote}

\begin{flushleft}
\end{flushleft}

\setcounter{equation}{0}
\section{Introduction}
Karatzas et al. \cite{KarLehShr87}, Cox and Huang \cite{CoxHua89} establish the static martingale method for the portfolio selection problem. This methodology identify the optimal
terminal wealth in closed form. Going one step further Ocone and Karatzas \cite{KAR2} represented the optimal terminal wealth using the Clark-Haussmann-Ocone formula. Detemple et al.
 \cite{DET1} provides new results on the structure of optimal portfolios including a portfolio decomposition formula in a more specialized model.

In some special cases people were able to derive portfolio decomposition formulas when the markets are not necessarily complete. Kim and Omberg \cite{KIM} look at an agent with CRRA preferences and a
 market model consisting of one risky asset $S(t)$ defined through
$$\frac{dS(t)}{S(t)}=\mu(t)dt+\sigma(t)dW(t),$$ and one risk-free asset with constant rate of return $r.$ The drift $\mu(t)$ and the volatility $\sigma(t)$ are diffusion processes.
The market price of risk $\theta(t)=\frac{\mu(t)-r}{\sigma(t)}$ follows an Ornstein-Uhlenbeck process. The analytical solution derived answers some qualitative questions: $1)$ when does the optimal portfolio hold more or less of the \textit{the myopic} component; $2)$ when does the optimal portfolio hold more of the risky asset at long horizon.

The purpose of this paper is to derive a new portfolio decomposition formula in complete markets. If $x$ is the agent's initial wealth, exogenously given,
 we show that if her wealth is topped up by a process $V_{x}(t)$ (which is described in Section $3$) then the optimal portfolio $\tilde{\pi}_{x}$ is robust with respect to stopping times (Theorem $3.1$), and we call it the \textit{the myopic portfolio}. If we fix a finite planning horizon $T$ and regard the random variable $V_{x}(T)$ as the payoff of a contingent claim, we can hedge it by trading in the financial market using some of the initial wealth. The \textit{hedging portfolio} $\bar{\pi}_{x}$ is obtained by Clark-Haussmann formula. An alternative way to find $\bar{\pi}_{x}$ is by means of Malliavin's calculus, i.e., the Karatzas-Ocone formula
$$V_{x}(T)=\mathbb{\tilde{E}}V_{x}(T)+\int_{0}^{T}\mathbb{\tilde{E}}\left[\left(D_{t}V_{x}(T)-V_{x}(T)\int_{t}^{T}D_{t}\tilde{\theta}(u)d\tilde{W}(u)\right)\bigg\vert\mathcal{F}(t)\right]d\tilde{W}(t). $$
However this formula works under the boundeness assumption on the market price of risk process $\tilde{\theta}$
 (see Theorem $2.5$ in \cite{KAR2}). In some models (see Kim and Omberg \cite{KIM}), $\tilde{\theta}$ is for instance an Ornstein-Uhlenbeck process which fails to remain bounded. 
 
 We choose another route and use the Clark-Haussmann formula, which gives $\bar{\pi}_{x}$ in terms of the Fr\'echet
 derivative of the functional $V_{x}(T)$ and the solution of a linearized SDE. It can be extended to apply to unbounded $\tilde{\theta}.$ Let us point out that this formula can be also employed to obtain optimal portfolios associated with option pricing.

The main mathematical result of this paper, Theorem $4.1,$ is a nontrivial extension of the Clark-Haussmann formula. We employ a judiciously chosen approximation-stopping procedure in order to represent $V_{x}(T).$ This in turn will give us \textit{the hedging portfolio}. The remainder of this paper is structured as follow. In Section $2$ we introduce the financial market model and describe the objective. Section $3$ derives \textit{the myopic portfolio}. Section $4$ deals with \textit{the hedging portfolio}. In order to accomplish this we extend Theorem $1$ from \cite{HA} to cover the functional $V_{x}(T).$ We conclude with an appendix containing some technical Lemmas.

\section{Model Description}
  
\subsection{Financial Market}
 We adopt a model for the financial market consisting of one bond and $d$ stocks. We work in discounted terms, that is the price
 of bond is constant and the stock price per share satisfy
 $$
dS_i(t)=S_i(t)\left[\alpha_i(t)\,dt
+\sum_{j=1}^n\sigma_{ij}(t)\,dW_j(t)\right],\quad0\leq
t\leq\infty,\quad i=1,\dots,n.
$$
Here $W=(W_{1},\cdots,W_{n})^{T}$ is a $n-$dimensional Brownian motion on a filtered probability space
$(\Omega,\{\sF_t\}_{0\leq t\leq T},\mathcal{F},\mathbb{P}),$ where $\{\sF_t\}_{0\leq t\leq T}$ is the completed filtration generated by $W.$ As usual $\{\alpha(t)\}_{t\in[0,\infty)}=\{(\alpha_i(t))_{i=1,\cdots,n}\}_{t\in[0,\infty)}$ is an $\mathbb{R}^{n}$ valued \textit{mean rate of return} process,
and\newline $\{\sigma(t)\}_{t\in[0,\infty)}=\{(\sigma_{ij}(t))_{i=1,\cdots,n}^{j=1,\cdots,n}\}_{t\in[0,\infty)}$ is an $n\times n-$matrix valued \textit{volatility} process,
 and are assumed progressively measurable with respect to $\{\sF_t\}_{0\leq t\leq T}.$  

\textit{Standing Assumption 2.1}
The matrix $\sigma(t)$ has full rank for every $t.$

This says that there are no redundant assets, and implies the existence of the inverse $(\sigma(t))^{-1}$ and the market price of risk process
\begin{equation}\label{1}\tilde{\theta}(t)=(\sigma(t))^{-1}\alpha(t). \end{equation}
All the processes encountered are defined on the fixed, finite interval $[0,T].$

\textit{Standing Assumption 2.2}
\begin{equation}\label{narde} \mathbb{E}\left[\exp\left(\frac{1}{2}\int_{0}^{T}\parallel\tilde{\theta}(u)\parallel^{2}\,du\right)\right]<\infty, \end{equation}
 where as usual $||\cdot||$ denotes the Euclidean norm in $\mathbb{R}^{n}.$ One can recognize this as the Novikov condition and it is sufficient to ensure that the
 stochastic exponential process 
 
  \begin{equation}\label{Nov}
  \tilde{Z}(t)=Z_{\tilde{\theta}}(t)\triangleq\exp\left\{-\int_{0}^{t}\tilde{\theta}^{T}(u)\,dW(u)-\frac{1}{2}\int_{0}^{t}\parallel\tilde{\theta}(u)\parallel^{2}\,du\right\}\end{equation}
is a (true) martingale. Moreover by the Girsanov theorem (Section $3.5$ in \cite{KAR}) \begin{equation}\label{Br}\tilde{W}(t)=W(t)+\int_{0}^{t}\tilde{\theta}(u)\,du \end{equation} is a Brownian motion under the equivalent martingale measure
  \begin{equation}\label{Q}
  \mathbb{\tilde{Q}}(A)\triangleq \mathbb{E}[\tilde{Z}(T){\bf{1}}_{A}],\qquad A\in \sF_T.\end{equation}
  Below we shall have occasion to write the process $\tilde{\theta}(t),$ cf. (\ref{1}), as a function of the process $W(t),$ i.e., $\tilde{\theta}(t,\omega)=\tilde{\Theta}(t,W(\cdot,\omega))$ a.s. 
Now define a mapping $\bar{\Theta}$ of $(C[0,T])^{n}$ into $(L^{\infty}[0,T])^{n}$ by $\bar{\Theta}(y)(t)\triangleq\tilde{\Theta}(t,y(\cdot)).$ Then $\bar{\Theta}$ is \textit{nonanticipative} in the sense that $\bar{\Theta}(y)(t)=\bar{\Theta}(z)(t)$
 for $y, z$ such that $y(s)=z(s)$ on $0\leq s\leq t$; this is equivalent to demanding $\tilde{\Theta}(t,W(\cdot))$ is $\{\sF_t\}$ -adapted. 
 
 \textit{Standing Assumption 2.3}
 \begin{equation}\label{as1}
 \sup_{t}\parallel\tilde{\Theta}(t,0)\parallel<\infty,
 \end{equation} 
and $\bar{\Theta}$ is Frech\'et differentiable with derivative $\bar{\Theta}'(y),$ i.e. for $y\in(C[0,T])^{n},$
 $\bar{\Theta}'(y)$ is a bounded linear operator mapping  $(C[0,T])^{n}$ into $(L^{\infty}[0,T])^{n}$ such that $$\parallel \bar{\Theta}(y+h)-\bar{\Theta}(y)-\bar{\Theta}'(y)h \parallel_{\infty}=o(\parallel h \parallel_{T}),$$  where $\parallel\cdot\parallel_{T}$ is the norm in $(C[0,T])^{n},$ i.e. $\parallel y\parallel_{T}=\sup_{0\leq t\leq T}\parallel y(t)\parallel,$ and $\parallel \cdot\parallel_{\infty}$ is the norm in $(L^{\infty}[0,T])^{n}.$ 
 
 The Riesz Representation Theorem gives, for fixed $t,$ the existence of a unique finite signed measure $\tilde{\mu}$  such that
  \begin{equation}\label{meas}[\bar{\Theta}'(y)h](t)=\int_{0}^{t}\tilde{\mu}(ds,y,t)h(s).\end{equation}

 \textit{Standing Assumption 2.4}
 We require that
 for some $\delta>0$ and constant $K_{\delta}$
 \begin{equation}\label{as2}
 |\bar{\Theta}'(y_1)-\bar{\Theta}'(y_2)|_{t}\leq
 K_{\delta}\parallel y_1-y_2 \parallel_{T}^{\delta}
 ,\end{equation}
 and
 
 \begin{equation}\label{as3}
  \sup_{y,t}|\bar{\Theta}'(y)|_{t}=\sup_{y,t}\mathrm{var}_{[0,t]}(\tilde{\mu}(\cdot,y,t))<\infty,  
 \end{equation}
 where
 $$|\bar{\Theta}'(y)|_{t}=\mathrm{var}_{[0,t]}(\tilde{\mu}(\cdot,y,t)). $$
 
 Recall that \begin{equation}\label{10o}dW(t)=d\tilde{W}(t)-\tilde{\Theta}(t,W(\cdot))dt=d\tilde{W}(t)-\bar{\Theta}(W)(t)dt,\end{equation} and by the above assumptions this SDE
 (where the unknown process is $W$) has a unique solution $W$ (see Theorem $6,$ page $249$ in \cite{Protter}). Hence $W$ and $\tilde{W}$ generate the same filtration
 $\{\sF_t\}_{0\leq t\leq T}.$
 
 In what follows we denote by $\mathbb{\tilde{E}}$ the expectation operator with respect to the probability measure $\mathbb{\tilde{Q}}.$ 
 
 \begin{lemma}\label{L:11}
For any nonnegative $\kappa$
\begin{equation}\label{qwe1a}
\mathbb{\tilde{E}}\parallel\bar{\Theta}(W)\parallel_{T}^\kappa=\mathbb{\tilde{E}}\sup_{t\leq T}\parallel \tilde{\Theta}(t,W(\cdot))\parallel^{\kappa}<\infty.
\end{equation} 
\end{lemma}
\noindent {\sc Proof:} See the Appendix. 
\begin{flushright}
$\square$
\end{flushright}

 One market model which fits into our framework is a stochastic volatility model, in which the market price of risk follows an Ornstein-Uhlenbeck process.
 The market consists of one bond and one stock whose price $S(t)$ is given by
 $$\frac{dS(t)}{S(t)} =\mu(t)\,dt
+\sigma(t)\,dW(t).$$ Let $U(t)$ be an Ornstein-Uhlenbeck process $$dU(t)=(\alpha-\beta U(t))dt+vdW(t),$$
In this model
 \begin{equation}\label{3z}\tilde{\theta}(t)=U(t)\end{equation} which is an unbounded process. Condition (\ref{as1}) is obviously satisfied and one can prove that (\ref{narde}) hold. 
In fact  \begin{equation}\label{e3}U(t)=e^{-\beta t}U(0)+\left(\frac{\alpha}{\beta}+\frac{v^{2}}{2\beta}\right)(1-e^{-\beta t})+ve^{-\beta t}\int_{0}^{t}e^{\beta u}dW(u),\end{equation}
whence

 $$ [\bar{\Theta}'(W)(\gamma)](t)  =v\left[\gamma(t)-\beta\int_{0}^{t}e^{\beta( u-t)}\gamma(u)\,du\,\right],$$ so (\ref{as2}) and (\ref{as3}) hold.

\subsection{Portfolio and wealth processes}

 A (self-financing) portfolio is defined as a pair $(x,\pi).$ The constant $x,$ exogenously given, is the initial value of the
 portfolio and $\pi=(\pi_{1},\cdots,\pi_{n})^{T}$ is a predictable $S-$ integrable process which specify how many units of the asset
 $i$ are held in the portfolio at time $t.$ The wealth process of such a portfolio is given by
 \begin{equation}\label{welath}
 X^{x,\pi}(t)=x+\int_{0}^{t}\pi(u)^{T}dS(u).
 \end{equation}

\subsection{Utility Function}
 A function $U:(0,\infty)\rightarrow \mathbb{R}$ strictly increasing and strictly concave is called a utility function.
 We restrict  ourselves to utility functions which are $4-$times continuous differentiable and satisfy the Inada conditions
 \begin{equation}\label{In}
 U'(0+)\triangleq\lim_{x\downarrow 0}U'(x)=\infty,\quad U'(\infty)\triangleq\lim_{x\uparrow{\infty}}U'(x)=0.
 \end{equation}
We shall denote by $I(\cdot)$ the (continuous, strictly decreasing) inverse of the marginal utility function $U'(\cdot),$ and by \eqref{In}
\begin{equation}\label{In1}
 I(0+)\triangleq\lim_{x\downarrow 0}I(x)=\infty,\quad I(\infty)\triangleq\lim_{x\uparrow{\infty}}I(x)=0.
 \end{equation}
\textit{Standing Assumption $2.1$}
\begin{equation}\label{0op}
y^{2}|I''(y)|\vee(-yI'(y))\vee I(y)<k_1y^{-\alpha}\quad\mbox{for\,\,every}\,\,y\in(0,\infty),
\end{equation}
for some $k_1>0,$ and $a\vee b=\max(a,b).$

\subsection{Objective}
 For a given initial positive wealth $x$ and a given utility function $U$ which satisfy the above assumptions, the optimal portfolio, $\hat\pi$, for (\ref{eq301}) is known to exist and can be obtained using the martingale representation theorem, cf. Karatzas and Shreve \cite{KarLehShr87}, 
\begin{equation}\label{eq301}
\sup_{\pi\in\mathcal{A}(x)}\mathbb{E}U(X^{x,\pi}(T))=\mathbb{E}U(X^{x,\hat{\pi}}(T)).
\end{equation} 
Here $\mathcal{A}(x)$ is the set of admissible portfolios. It is defined by
\begin{equation}\label{admis} 
\mathcal{A}(x)\triangleq\left\{\pi|
X^{x,\pi}(t)>0,\,0\leq t\leq T,\,\,\,\mathbb{E}[U(X^{x,\pi}(T)]^{-}<\infty\right\},
\end{equation}
where $a^{-}\triangleq\max\{-a,0\}.$ Of course the agent may be uncertain about her investment horizon so she would like a robust or time-consistent optimal policy, i.e., for any stopping time $\tau\leq T$
\begin{equation}\label{Probust}
\tilde{\pi}\in{\arg\max_{\pi}}\,\mathbb{E}U(X^{x,\pi}(\tau)).
\end{equation}
 Put differently, two agents with the same preferences, living on $[0,\tau]$ and $[\tau,T]$ use the same optimal portfolio, $\tilde{\pi}$, as one agent living on $[0,T].$ The problem can also be viewed in terms of the consistency problem described in \cite{EKE}.  
  
A solution, $\tilde\pi$, of (\ref{Probust}) does not exist unless the utility is logarithmic. However, as shown in \cite{PIR}, if we top up the wealth by a finite variation process, then there exists $\tilde{\pi}$ such that for any stopping time $\tau\leq T$,
\[
\tilde{\pi}={\arg\max_{\pi}}\,\mathbb{E}U(X^{x,\pi}(\tau)+V_{x}(\tau)),
\] 
where $V_{x}$ is a finite variation process which depends on the utility function. 
 We regard this process as a measure of ``time-inconsistency'' of the investment policy of an
 agent due to her non-logarithmic utility. Indeed for log, $V_{x}=0,$ and $\hat{\pi}=\tilde{\pi}.$ The process  $V_{x}$ can also be seen as a ``risk measure of time consistency" because at any time is it the amount of money needed
 to be added to the investors's wealth to yield time consistency of the investor's optimal wealth. 
Our aim then is to decompose $\hat\pi$ into the process $\tilde{\pi}$ and a hedging component, $\bar{\pi}$ (which depends on the stopping time $T$ here), such that
\begin{equation}\label{dec}
\hat{\pi}=\tilde{\pi}+\bar{\pi}.
\end{equation}
This requires a division of the initial wealth into a part which is invested according to $\tilde\pi$ and the rest which is invested according to $\bar\pi$ to generate the corresponding $V_x(T)$. In case of non-uniqueness, we invest as much as possible in the myopic part, i.e. we minimize the time inconsistency of the optimal portfolio.    

In comparison to other portfolio decomposition formulas, ours shows the robustness of the myopic component $\tilde{\pi}$ with respect to all the stopping times. Moreover the structure of the optimal portfolio reflects the time inconsistency due to the agent's utility function through the hedging component $\bar{\pi}.$ We find such $\tilde\pi$ in Section 3 and observe that it is myopic, i.e. does not depend on the future evolution of stock prices.

\section{The Myopic Portfolio}

In this section we find the myopic portfolio $\tilde{\pi}$. In order to accomplish this we solve the problem of maximizing expected utility of
the final wealth adjusted by a finite variation process. This makes $\tilde{\pi}$ robust with respect to all stopping times. Indeed let us consider
\begin{equation}\label{q3}
V_{x}(t)=\int_{0}^{t}F(U'(x)\tilde{Z}(u))||\tilde{\theta}(u)||^{2}\,du,\quad 0\leq t\leq T,
\end{equation}
where
\begin{equation}\label{0q3}
F(z)=\frac{1}{2}I''(z)z^{2}+I'(z)z.
\end{equation} 
and the corresponding set of admissible portfolios 
\begin{equation}\label{leq4}
\mathcal{A}_{V}(x)\triangleq\{\pi|\,\,X^{x,\pi}(t)+V_{x}(t)>0,\,\,0\leq t\leq T\}.
\end{equation}

\begin{theorem}\label{main1}
Let $\tau$ be a stopping time. The optimal portfolio process for maximizing expected utility of
the final wealth adjusted by $V_{x},$ i.e.,
\begin{equation}\label{eq3}
{\sup_{\pi\in\mathcal{A}_{V}(x)}}\mathbb{E}U(X^{x,\pi}(\tau)+V_{x}(\tau))=\mathbb{E}U(X^{x,\tilde{\pi}}(\tau)+V_{x}(\tau)),
\end{equation} 
is given by 
\begin{equation}\label{q34}
(\tilde{\pi}_{x}(t))_{i}=-\frac{1}{S_{i}(t)}((\sigma^{T}(t))^{-1}U'(x)I'(U'(x)\tilde{Z}(t))\tilde{Z}(t)\tilde{\theta}(t))_{i}, \quad 0\leq t\leq T,\,\, i=1,\cdots,d. 
\end{equation}
\end{theorem} 

\noindent {\sc Proof:} See the Appendix.

\begin{flushright}
$\square$
\end{flushright}

\begin{remark}
If the utility is logarithmic, i.e.,
$U(x)=\log{x},$ then $V_{x}(t)\equiv 0$. Being optimal for the logarithmic utility, the vector $\tilde{\pi}_{x}$ satisfies
\begin{equation}\label{0i}(\tilde{\pi}_{x}(t))_{i}=\frac{(\zeta_{M}(t))_{i}X^{\tilde{\pi}}(t)}{S_{i}(t)},\quad i=1,\dots,n,\end{equation} with $\zeta_{M}(t)\triangleq(\sigma(t)\sigma^\mathrm{T}(t))^{-1}\alpha(t)$ the Merton proportion. The future evolution of $S$ does not enter in the formula (\ref{q34}) and (\ref{0i}), hence we refer to $\tilde{\pi}_{x}$  as \textit{the myopic component}.
\end{remark}

\section{The Hedging Portfolio}
The correction term process
\begin{equation}\label{3q3}
V_{x}(T)=\int_{0}^{T}F(U'(x)\tilde{Z}(u))||\tilde{\theta}(u)||^{2}\,du,\quad 0\leq t\leq T,
\end{equation}
can be viewed as a contingent claim and requires some of the initial wealth to be hedged. It admits the following representation
\begin{equation}\label{2q3}
V_{x}(T)=\int_{0}^{T}\bar{\pi}_{x}(u)^{T}\,dS(u)+\tilde{\mathbb{E}}V_{x}(T),
\end{equation} 
for some process $\bar{\pi}_{x}.$ In the reminder of this section we show how to compute $\bar{\pi}_{x}$ explicitly.

The martingale ${\bar{V}_{x}}(t)\triangleq\tilde{\mathbb{E}}[V_{x}(T)|\sF_{t}],$\,\,admits the stochastic integral representation
 \begin{equation}\label{8i}
\bar{V}_{x}(t)=\tilde{\mathbb{E}}V_{x}(T)+\int_{0}^{t}\beta^{T}(u)d\tilde{W}(u),\qquad 0\leq t\leq T,
\end{equation}
 for some $\sF_t-$adapted process $\beta(\cdot)$ which satisfies $\int_{0}^{T} ||\beta(u)||^{2}\,du<\infty$ a.s. (e.g., \cite{KAR1}, Lemma $1.6.7$).

In the light of this we want $$\int\beta^{T}\,d\tilde{W}=\int\bar{\pi}^{T}_{x}\,dS.$$ Let $A$ be the $n\times n$ matrix with the entries $A_{ij}=S_i\sigma_{ij}.$ Since the volatility matrix $\sigma$ has linearly independent rows, the matrix $A$ has linearly independent rows, i.e., $\mbox{rank}\,A=n.$ Therefore
\begin{equation}\label{rep}
\bar{\pi}_{x}= (A^{T})^{-1}\beta
\end{equation}

Let us notice that the representation formula in (\ref{8i})
takes place under the probability measure $\tilde{\mathbb{Q}}.$ However by Theorem $14$ page $60$ in \cite{Protter} this takes place under $\mathbb{P},$ since $\mathbb{P}\sim \tilde{\mathbb{Q}}.$  The process $\beta(t)$ of \eqref{8i} can be computed explicitly by Haussmann's formula. Let us define the functional $L:C[0,T]\times(C[0,T])^{n}\longrightarrow\mathbb{R}$ as
  
   \begin{equation}\label{eeu}
   L(z,w)\triangleq\int_{0}^{T}F(U'(x)z(u))||\bar{\Theta}(w)(u)||^{2}\,du,
   \end{equation}
   with $\bar{\Theta}(W)(u)=\tilde{\theta}(u,\omega).$ Then $L(\tilde{Z},W)={V}_{x}(T),$ and it can be shown that
   \begin{eqnarray}\label{pis}
   L'(z,w)(v_1,v_2)&=&\int_{0}^{T}\big\{[F'(U'(x)z(u))]\parallel\bar{\Theta}(w)(u)\parallel^{2}v_1(u)\\\notag& & \mbox{}+2F(U'(x)z(u))\bar{\Theta}^{T}(w)(u)[\bar{\Theta}'(w)v_2](u)\big\}du \\\notag&=&\int_{0}^{T}\mu(dt,z,w)(v_1(t),v_2(t)^{T})^{T},\,\,\,\mbox{a.s.,}\notag
   \end{eqnarray}
   for some measure $\mu.$ Let $Y=(\tilde{Z},W);$ thus
 $$Y(t)=Y(0)+\int_0^t f(u,Y(\cdot))\,du+\int_0^t g(u,Y(\cdot))d\tilde{W}(u),\quad 0\leq t\leq T. $$
  Here with $y=(y_1,y_2),$ $y_1$ a scalar process and $y_2$ a $n$-dimensional process, 
 $$f(u,y)=\begin{pmatrix} y_1\parallel\tilde{\Theta}(u,y_2)\parallel^2 \\-\tilde{\Theta}(u,y_2) \end{pmatrix} $$
 $$g(u,y)=\begin{pmatrix} -y_1\tilde{\Theta}^{T}(u,y_2)\\I_n \end{pmatrix}, $$
 where $I_n$ is the $n\times n$ identity matrix. Let $\Phi(t,s)$ be the unique solution of the linearized equation 
$$d\Phi(t,s)=\left[\frac{\partial f}{\partial y}(t,Y(\cdot))\Phi(t,s)\right](t)dt+\left[\frac{\partial g}{\partial y}(t,Y(\cdot))\Phi(t,s)\right](t)d\tilde{W}(t),\qquad t>s,$$
   $$\Phi(s,s)=I_{n+1},\,\,\mbox{and}\,\,\,\,\Phi(t,s)=O_{n+1}\,\,\,\,\mbox{for}\,\,\,\,0\leq t<s.$$
     At this point one may wonder about the existence and uniqueness of $\Phi.$ The matrix  process $\Phi$ has $j-$th column $\Phi^{j}=(\Phi^{1,j} ,(\Phi^{2,j})^{T})^{T}.$
The scalar process $\Phi^{1,j}$ satisfies
\begin{eqnarray*}
d\Phi^{1,j}(t,s)&=&\left(\parallel\bar{\Theta}(W)(t)\parallel^2\Phi^{1,j}(t,s)+2\tilde{Z}(t)\bar{\Theta}^{T}(W)(t)[\bar{\Theta}'(W)\Phi^{2,j}(\cdot,s)](t)\right)dt\\ &  & \mbox{} - \left(\bar{\Theta}^{T}(W)(t)\Phi^{1,j}(t,s)+\tilde{Z}(t)[\tilde{\Theta}'(W)\Phi^{2,j}(\cdot,s)](t)\right)d\tilde{W}(t).
\end{eqnarray*}
The $n-$dimensional vector process $\Phi^{2,j}$ satisfies
\begin{equation}\label{5}
\frac{d\Phi^{2,j}(t,s)}{dt}=-\left[\bar{\Theta}'(W)\Phi^{2,j}(\cdot,s)\right](t)=-\int_{s}^{t}\tilde{\mu}(du,W,t)\Phi^{2,j}(u,s),
\end{equation}
where $\tilde{\mu}$ is the measure defined in \eqref{meas}. Let us notice that $\Phi^{2,1}\equiv0.$  For $j>1$ and some constant $K,$
\begin{eqnarray*}
\parallel\Phi^{2,j}(t,s)-e_{j-1} \parallel&=&\parallel\Phi^{2,j}(t,s)-\Phi^{2,j}(s,s) \parallel\\&=&\parallel\int_{s}^{t}-[\bar{\Theta}'(W)\Phi^{2,j}(\cdot,s)] (u)\,du\parallel\\&\leq&\int_{s}^{t}K\sup_{s\leq v\leq u}\parallel\Phi^{2,j}(v,s)\parallel\,du 
,\end{eqnarray*} where  $e_{i}$ is the $i^{th}$ column of $I_{n}$ and the last inequality comes from (\ref{as3}). Therefore
$$ \sup_{s\leq v\leq t}\parallel\Phi^{2,j}(v,s)\parallel\leq 1+K\int_{s}^{t}\sup_{s\leq v\leq u}\parallel\Phi^{2,j}(v,s)\parallel\,du ,$$ hence by Gronwall's inequality
\begin{equation}\label{oop}
 \sup_{s\leq v\leq t}\parallel\Phi^{2,j}(v,s)\parallel\leq e^{K(t-s)}.
\end{equation}
Existence and uniqueness of the process $\Phi$ it is now straightforward. Let us define
 \begin{equation}\label{12}\lambda(t,\omega)\triangleq\left[\int_{t}^{T}\mu(du,\tilde{Z}(\omega),W(\omega))\Phi(u,t,\omega)\right]g(t,(\tilde{Z}(\omega),W(\omega))), \end{equation}
 with $\mu$ the measure of \eqref{pis}.

 \begin{theorem}\label{main22}
  (Clark-Haussmann formula)
 \begin{equation}\label{23}L(\tilde{Z}(\cdot),W(\cdot))=\int_{0}^{T}\mathbb{\tilde{E}}(\lambda(t)\vert\mathcal{F}_{t})d\tilde{W}(t)+\mathbb{\tilde{E}}L[\tilde{Z}(\cdot),W(\cdot)]. \end{equation} \end{theorem}

\noindent {\sc Proof:} See the Appendix.
\begin{flushright}
$\square$
\end{flushright}

  The process $\mathbb{\tilde{E}}(\lambda(t)\vert\mathcal{F}_{t})$ satisfies the integrability condition (see (\ref{show2})) and the discussion following it) to grant
$$\mathbb{\tilde{E}}(\lambda(t)\vert\mathcal{F}_{t})=\beta^{T}(t),$$
with $\beta$ of \eqref{8i}.

\textit{Standing Assumption $5.1$} 
\begin{equation}\label{bg}
-\infty<x_{*}\triangleq\inf\{z\geq0|z+\tilde{\mathbb{E}}V_{z}(T)=x\}<\infty.
\end{equation}

\begin{remark}
Let us compute the process $V_{x}$ for different utility functions using \eqref{q3}. 
In the case of an exponential utility, i.e., $U(x)=-e^{-ax}$:
$$V_{x}(t)=\frac{1}{a}\int_{0}^{t}||\tilde{\theta}(u)||^{2}\,du. $$
If the utility is CRRA, i.e., $U(x)=\frac{x^{p}}{p}$:
$$V_{x}(t)=\frac{xp}{2(p-1)^{2}}\int_{0}^{t}[\tilde{Z}(u)]^{\frac{1}{p-1}}||\tilde{\theta}(u)||^{2}\,du.$$
Thus for exponential utility assumption \eqref{bg} holds true. As for power utility we only need the weaker assumption
$$1+\frac{p}{2(p-1)^{2}}\mathbb{\tilde{E}}\int_{0}^{T}[\tilde{Z}(u)]^{\frac{1}{p-1}}||\tilde{\theta}(u)||^{2}\,du\neq 0.$$
\end{remark}

\begin{theorem}\label{main3}
Starting with the initial wealth $x,$ the agent should invest $\tilde{\pi}_{x_{*}}+\bar{\pi}_{x_{*}}$ in stocks, where
$$(\tilde{\pi}_{x_{*}}(t))_{i}=-\frac{1}{S_{i}(t)}((\sigma^{T}(t))^{-1}U'(x_{*})I'(U'(x_{*})\tilde{Z}(t))\tilde{Z}(t)\tilde{\theta}(t)^{T})_{i}, \quad 0\leq t\leq T,\,\, i=1,\cdots,d.    $$
is \textit{the myopic portfolio}, $\bar{\pi}_{x_{*}}$ of \eqref{rep} is \textit{the hedging portfolio}, and $x_{*}$ is given by \eqref{bg}. This investment strategy is optimal, i.e.,
\begin{equation}\label{eq0}
{\sup_{\pi\in\mathcal{A}(x)}}\mathbb{E}U(X^{x,\pi}(T))=\mathbb{E}U(X^{x,\tilde{\pi}+\bar{\pi}}(T)).
\end{equation}
\end{theorem}

\noindent {\sc Proof:} Since
$$X^{x,\tilde{\pi}+\bar{\pi}}(T)=I(U'(x_{*})\tilde{Z}(T)),$$
one can argue as in the proof of Theorem \ref{main1} to obtain optimality.

\begin{flushright}
$\square$
\end{flushright}

\section{Appendix}

{\bf{Proof of Lemma \ref{L:11}}}: Since
 $$\bar{\Theta}(y)(t)=\bar{\Theta}(0)(t)+\int_{0}^{1}[\bar{\Theta}'(uy)y](t)\,du ,$$ (\ref{as1}) and (\ref{as3}) yield a constant $K_1$ such that
 \begin{eqnarray}\label{ver}
\| [\bar{\Theta}'(uy)y](t) \|&\leq& K_1+\left\|\int_{0}^{1}\int_{0}^{t}\tilde{\mu}(ds,y,uy)y(s)\,du\right\|
\\\notag&\leq& K_1+\int_{0}^{1}\mathrm{var}_{[0,t]}(\tilde{\mu}(\cdot,y,uy)\|y\|_{t}\,du\\&\leq& K_1(1+\|y\|_{t}).\notag
\end{eqnarray}
 Moreover
 \begin{eqnarray*}\parallel W(t)\parallel&\leq& \parallel\tilde{W}(t)\parallel+\int_{0}^{t}\parallel\bar{\Theta}(W)(u)\parallel\,du
\\&\leq&\parallel\tilde{W}(t)\parallel+\int_{0}^{t}K_1(1+\parallel W \parallel_{u})\,du.
\end{eqnarray*}
Therefore \begin{equation}\label{poi}\parallel W\parallel_{T}^{\kappa}\leq K_2+K_3(\parallel\tilde{W}\parallel_{T}^{\kappa}+\int_{0}^{T}\parallel W\parallel_{t}^{\kappa}\,dt),\end{equation} for some constants $K_2$ and $K_3.$ The process $\tilde{W}$ being a Brownian motion under $\tilde{Q},$ has finite moments, hence (\ref{poi}) and Gronwall's inequality prove (\ref{qwe1a}). 
  \begin{flushright}
$\diamond$
\end{flushright}

\noindent {\bf{Proof of Theorem \ref{main1}}}:
For $\pi\in\mathcal{A}_{V}(x),$ the process $X^{x,\pi}$ is a continuous $\mathbb{\tilde{Q}}$ local martingale. Moreover is a supermartingale, since
$$\mathbb{\tilde{E}}[\sup_{0\leq t\leq T}|V_{x}(t)|]<\infty, $$ 
which is a consequence of \eqref{0op} and \eqref{qwe1a}. Problem $3.26,$ p. $20$ in \cite{KAR} yields  
\begin{equation}\label{i9}
\mathbb{E}[\tilde{Z}(\tau)X^{x,\pi}(\tau)]\leq x.
\end{equation} 
It turns out that with $\tilde{\pi}$ of \eqref{q34}
$$X^{x,\tilde{\pi}}(t)=x-\int_{0}^{t}U'(x)I(U'(x)\tilde{Z}(u))\tilde{Z}(u)\tilde{\theta}^{T}(u)d\tilde{W}(u).$$
 The assumption \eqref{0op} in conjunction with \eqref{qwe1a} make the process $X^{x,\tilde{\pi}}(t)$
a (true) martingale. Hence for any stopping time $\tau\leq T$ 
\begin{equation}\label{1i9}
\mathbb{E}[\tilde{Z}(\tau)X^{x,\tilde{\pi}}(\tau)]=x,
\end{equation}  
by Problem $3.26,$ p. $20$ in \cite{KAR}.
Direct computations show
\begin{equation}\label{eu1}
X^{x,\tilde{\pi}}(\tau)+V_{x}(\tau)=I(U'(x)\tilde{Z}(\tau)),
\end{equation}
with $V_{x}$ of \eqref{q3}. This and the concavity of the utility function yield

  \begin{eqnarray*}
  U(X^{x,\pi}(\tau)+V_{x}(\tau))-U(X^{x,\tilde{\pi}}(\tau)+V_{x}(\tau))&\leq& U'(X^{x,\tilde{\pi}}(\tau)+V_{x}(\tau))(X^{x,\pi}(\tau)-X^{x,\tilde{\pi}}(\tau))\\&=&U'(x)\tilde{Z}(\tau)(X^{x,\pi}(\tau)-X^{x,\tilde{\pi}}(\tau)).  \end{eqnarray*}
Taking expectation and using \eqref{i9} and \eqref{1i9} we conclude. 
\begin{flushright}
$\diamond$
\end{flushright}

\noindent {\bf{Proof of Theorem \ref{main22}}}: Let us notice that we cannot apply Haussmann's formula right away because hypothesis $H_3$ in \cite{HA} fails.
 However we can get around this by an approximation argument. Let $\phi_k$ be a sequence of bounded differentiable
 functions on $\mathbb{R}$ with H\"older continuous (of order $\delta$  (see (\ref{as2}))) derivatives, such that $\phi_k(x)=x,$ if $|x|\leq k,$\,\,  $|\phi_k(x)|\leq |x|,$\,\, $|\phi'_k(x)|\leq1,$ and define
 $$\phi_k(\tilde{\theta})\triangleq (\phi_k(\tilde{\theta_{1}}),\cdots,\phi_k(\tilde{\theta_{n}})),$$
 
 $$f_{k}(u,y)=\begin{pmatrix}\phi_k(y_1)\parallel\phi_{k}(\tilde{\theta}(u,y_2))\parallel^2\\-\tilde{\theta}(u,y_2) \end{pmatrix},\qquad g_{k}(u,y)=\begin{pmatrix}-\phi_k(y_1)\phi_k(\tilde{\theta}(u,y_2))\\I_n\end{pmatrix}. $$

 It is easily seen that for fixed $k,$ the functions $f_k(\cdot,\cdot)$ and $g_k(\cdot,\cdot)$ satisfy
 the conditions of Theorem\,$1$ in \cite{HA}. Let us denote by $Y_{k}(t)$ the unique strong solution of the SDE 
 $$dY(t)=f_{k}(t,Y(\cdot))dt+g_{k}(t,Y(\cdot))d\tilde{W}(t),$$ (the existence of a unique strong solution follows since
  $f_k(\cdot,\cdot)$ and $\sigma_k(\cdot,\cdot)$ are Lipschitz). In fact $Y_{k}(t)=(\tilde{Z}_{k}(t),W(t))$ for some process $\tilde{Z}_{k}(\cdot).$ It turns out that the process  $\tilde{Z}_{k}(\cdot)$ is strictly positive.
  Indeed when $\tilde{Z}_{k}(\cdot)$ gets close to zero it satisfies 
   $$d\tilde{Z}_{k}=\tilde{Z}_{k}(\parallel\phi_{k}(\tilde{\theta})\parallel^2\,dt-\phi_{k}(\tilde{\theta})d\tilde{W}(t)) ,$$ hence the positivity.
As in \cite{HA} we need $L:(C[0,T])^{n+1}\longrightarrow\mathbb{R}$ and
   \begin{equation}\label{sid} |L'(y_1)-L'(y_2)|\leq K(1+\parallel y_1\parallel_{T}^{\beta})(1+\parallel y_2 \parallel_{T}^{\beta})\parallel y_1-y_2\parallel_{T}^{\rho},\end{equation} for some positive $K, \beta, \rho.$  The assumption (\ref{as2}) yields:
   $$\parallel\bar{\Theta}(y_1)-\bar{\Theta}(y_2)\parallel_{T}\leq K\parallel y_1-y_2\parallel_{T}+K_{\delta}(\parallel y_1\parallel_{T}\wedge\parallel y_2\parallel_{T})\parallel y_1-y_2\parallel^{\delta}_{T},$$ for some constant $K.$
   This together with (\ref{as2}), (\ref{as3}),(\ref{0op}) and (\ref{ver}) prove (\ref{sid}).

   Following \cite{HA}, let $\Phi_k(t,s)$ be the $(n+1)\times (n+1)$ matrix which solves the linearized equation $$dZ(t)=\left[\frac{\partial f_{k}}{\partial x}(t,Y_{k}(\cdot))Z(\cdot)\right](t)dt+\left[\frac{\partial g_k}{\partial x}(t,Y_{k}(\cdot))Z(\cdot)\right](t)d\tilde{W}(t),\qquad t>s,$$
 $$\Phi_{k}(s,s,\omega)=I_{n+1},\,\,\mbox{and}\,\,\,\,\Phi_{k}(t,s,\omega)=O_{n+1}\,\,\,\,\mbox{for}\,\,\,\,0\leq t<s,$$ with $O_{n+1}$ the $(n+1)\times(n+1)$ matrix with zero entries.
 Next with $\mu$ of \eqref{pis} we define $$\lambda_{k}(t,\omega)\triangleq\left[\int_{t}^{T}\mu(du,(\tilde{Z}_{k}(\omega),W(\omega)))\Phi_k(u,t,\omega)\right]g_{k}(t,(\tilde{Z}_{k}(\omega),W(\omega))). $$ Theorem $1$ in \cite{HA} gives the following representation (the Clark-Haussmann formula)
\begin{equation}\label{09}L(\tilde{Z}_{k}(\cdot),W(\cdot))=\int_{0}^{T}\mathbb{\tilde{E}}(\lambda_{k}(t)\vert\mathcal{F}_{t})d\tilde{W}(t)+\mathbb{\tilde{E}}L[\tilde{Z}_{k}(\cdot),W(\cdot)]\end{equation}
 
  At this point we need some auxiliary Lemmas.
  \begin{lemma}\label{L:0}
For every real number $r$
\begin{equation}\label{ris}
\mathbb{\tilde{E}}\sup_{0\leq t\leq T}\tilde{Z}_{k}^{r}(t)<\infty,
\end{equation}
uniformly in $k$ and
\begin{equation}\label{2}
\mathbb{\tilde{E}}\sup_{0\leq t\leq T}\tilde{Z}^{r}(t)<\infty.
\end{equation}
\end{lemma}
\noindent {\sc  Proof:} 
\vspace{0.5cm}
 We prove \ref{ris} and \ref{2} follows similarly. If $r\neq0$ is a real number, It\^o's Lemma gives
 \begin{eqnarray*}\tilde{Z}_{k}^{r}(t)&=&1+\int_{0}^{t}r\tilde{Z}_{k}^{r-1}(u)\phi_k(\tilde{Z}_{k}(u))\parallel\phi_{k}(\bar{\Theta}(W)(u))\parallel^2\,du\\& & \mbox{}+\int_{0}^{t}r(r-1)\tilde{Z}_{k}^{r-2}(u)\phi_k^{2}(\tilde{Z}_{k}(u))\parallel\phi_{k}(\bar{\Theta}(W)(u))\parallel^2\,du\\& & \mbox{}+\int_{0}^{t}r\tilde{Z}_{k}^{r-1}(u)\phi_k(\tilde{Z}_{k}(u))\phi_{k}(\bar{\Theta}(W)(u))^{T}\,d\tilde{W}(u).\end{eqnarray*}
 Consequently
\begin{eqnarray*}\tilde{Z}_{k}^{4r}(t)&\leq& K_{1}\bigg[1+\left|\int_{0}^{t}r\tilde{Z}_{k}^{r-1}(u)\parallel\phi_k(\bar{\Theta}(W)(u))\parallel^2\,du\right|^{4}\\ & & \mbox{}+ \left|\int_{0}^{t}r(r-1)\tilde{Z}_{k}^{r-2}(u)\phi_k^{2}(\tilde{Z}_{k}(u))\parallel\phi_{k}(\bar{\Theta}(W)(u))\parallel^2\,du\right|^{4}\\ & & \mbox{}+\left|\int_{0}^{t}r\tilde{Z}_{k}^{r-1}(u)\phi_k(\tilde{Z}_{k}(u))\phi_{k}(\bar{\Theta}(W)(u))^{T}\,d\tilde{W}(u)\right|^{4}\bigg].\end{eqnarray*}
 Due to  $|\phi_k(x)|\leq |x|,$ H\"older's inequality and (\ref{qwe1a}) imply 
\begin{eqnarray*}
\lefteqn{\mathbb{\tilde{E}}\left|\int_{0}^{t}\tilde{Z}_{k}^{r-1}(u)\phi_k(\tilde{Z}_{k}(u))\parallel\phi_{k}(\bar{\Theta}(W)(u))\parallel^2\,du\right|^{4}}\\ & & \mbox{} \leq   \mathbb{\tilde{E}}\left|\int_{0}^{t}\tilde{Z}_{k}^{r}(u)\parallel \bar{\Theta}(W)(u)\parallel^2\,du\right|^{4}\\& & \mbox{} \leq  \left(\mathbb{\tilde{E}}\int_{0}^{t}\tilde{Z}_{k}^{4r}(u)\,du\right)\left(\mathbb{\tilde{E}}\int_{0}^{t}\parallel \bar{\Theta}(W)(u)\parallel^{\frac{8}{3}}\,du \right)^{3}\\& & \mbox{} \leq  K_{2}\,\mathbb{\tilde{E}}\int_{0}^{t}\sup_{0\leq v\leq u}\tilde{Z}^{4r}_{k}(v)\,du.
\end{eqnarray*}
 Similarly
 $$ \mathbb{\tilde{E}}\left|\int_{0}^{t}\tilde{Z}_{k}^{r-2}(u)\phi_k^{2}(\tilde{Z}_{k}(u))\parallel\phi_{k}(\bar{\Theta}(W)(u))\parallel^2\,du\right|^{4} \leq K_{2}\,\mathbb{\tilde{E}}\int_{0}^{t}\sup_{0\leq v\leq u}\tilde{Z}^{4r}_{k}(v)\,du.$$

The Burkholder-Davis-Gundy inequality implies
\begin{eqnarray*}
\lefteqn{\mathbb{\tilde{E}}\sup_{0\leq t\leq T}\left|\int_{0}^{t}\tilde{Z}_{k}^{r-1}(u)\phi_k(\tilde{Z}_{k}(u))\phi_{k}(\bar{\Theta}(W)(u))^{T}\,d\tilde{W}(u)\right|^{4}}\\& & \mbox{} \leq \mathbb{\tilde{E}}\left(\int_{0}^{t}\tilde{Z}_{k}^{2(r-1)}(u)\phi_k(\tilde{Z}_{k}(u))^{2}\parallel\phi_{k}(\bar{\Theta}(W)(u))\parallel^{2}\,du\right)^{2}.
\end{eqnarray*}
Moreover, $|\phi_k(x)|\leq |x|,$ H\"older's inequality and (\ref{qwe1a}) give
\begin{eqnarray*}
\lefteqn{\mathbb{\tilde{E}}\left(\int_{0}^{t}\tilde{Z}_{k}^{2(r-1)}(u)\phi_k(\tilde{Z}_{k}(u))^{2}\parallel\phi_{k}(\bar{\Theta}(W)(u))\parallel^{2}\,du\right)^{2}}\\& &\mbox{}\leq\mathbb{\tilde{E}}\left(\int_{0}^{t}\tilde{Z}_{k}^{2r}(u)\parallel(\bar{\Theta}(W)(u))\parallel^{2}\,du\right)^{2}\\& &\mbox{}\leq\left(\mathbb{\tilde{E}}\int_{0}^{t}\tilde{Z}_{k}^{4r}(u)\,du\right)\left(\mathbb{\tilde{E}}\int_{0}^{t}\parallel\bar{\Theta}(W)(u)\parallel^{4}\,du\right)\\& &\mbox{}\leq K_{3}\,\mathbb{\tilde{E}}\int_{0}^{t}\sup_{0\leq v\leq u}\tilde{Z}^{4r}_{k}(v)\,du.\end{eqnarray*}
The above arguments show
$$\mathbb{\tilde{E}}\sup_{0\leq t\leq T}\tilde{Z}_{k}^{4r}(t)\leq K_{4}\int_{0}^{T}\mathbb{\tilde{E}}\sup_{0\leq v\leq u}\tilde{Z}_{k}^{4r}(v)\,du.$$
Finally Gronwall's inequality proves (\ref{ris}).
\begin{flushright}
$\square$
\end{flushright}

With the notations $a\vee b\triangleq\max(a,b)$
 and $|\tilde{\theta}|\triangleq\max\{|\tilde{\theta}_{i}|:\,i=1,\cdots,n\},$ let us define the following sequence of stopping times
  $$\tau_{k}\triangleq\inf\{s\leq T,\,\,\mbox{such\,that},\,\,\tilde{Z}_{k}(s)\vee|\bar{\Theta}(W)(s)|\geq k\} $$
 \begin{lemma}\label{L:1}
 $\tilde{Z}_{k}(t)=\tilde{Z}(t)$ for $t\leq\tau_{k}.$ Moreover $\tau_{k}\uparrow T$\,\,$\mathbb{P}$\,a.s. 
\end{lemma}

\noindent {\sc  Proof:} 
\vspace{0.5cm}
  Let us notice that on $[0,\tau_{k}],$\, $f_k=f$ and $g_k=g,$ hence $\tilde{Z}_{k}(t)=\tilde{Z}(t)$ (since it satisfies the same SDE ). Therefore \begin{eqnarray}\tau_{k}&\triangleq&\inf\{s\leq T,\,\,\mbox{such\,that},\,\,\tilde{Z}_{k}(s)\vee|\bar{\Theta}(W)(s)|\geq k\}\\&=&\inf\{s\leq T,\,\,\mbox{such\,that},\,\,\tilde{Z}(s)\vee|\bar{\Theta}(W)(s)|\geq k\}.\end{eqnarray} Thus for $t\leq\tau_k,$ one has $\tilde{Z}(t)\vee|\bar{\Theta}(W)(t)|\leq k,$ so $t\leq\tau_{k+1}.$
   This proves $\tau_{k+1}\geq\tau_{k}.$ For a fixed path $\omega,$ by (\ref{2}) and (\ref{qwe1a}) ${\sup_{0\leq s\leq T}}\tilde{Z}(s)\vee|\bar{\Theta}(W)(s)|\leq K(\omega),$ showing that  $\tau_{k}\uparrow T$\,\,$\mathbb{P}$\,\,a.s. 
\begin{flushright}
$\square$
\end{flushright}
\vspace{0.5cm}

\begin{lemma}\label{L:3}
 $$\mathbb{\tilde{E}}L(\tilde{Z}_{k}(\cdot),W(\cdot)) \longrightarrow\mathbb{\tilde{E}}L(
\tilde{Z}(\cdot),W(\cdot)),$$
as $k\rightarrow\infty,$ and by passing to a subsequence
$$L(\tilde{Z}_{k_{n}}(\cdot),W(\cdot))\longrightarrow L(\tilde{Z}(\cdot),W(\cdot))\quad \mbox{a.s.}$$
\end{lemma}
\noindent {\sc   Proof:} 
\vspace{0.5cm}

$$\mathbb{\tilde{E}}L(\tilde{Z}_{k}(\cdot),W(\cdot)) \longrightarrow\mathbb{\tilde{E}}L(\tilde{Z}(\cdot),W(\cdot)),$$ if and only if (cf. \eqref{eeu})
$$\mathbb{\tilde{E}}\int_{0}^{T}\left[F(U'(x)\tilde{Z}_{k}(u))\parallel \bar{\Theta}(W)(u) \parallel^2-   F(U'(x)\tilde{Z}(u))\parallel\bar{\Theta}(W)(u) \parallel^2\right]\,du\longrightarrow 0.$$ Let us notice that for $u<T$
\begin{eqnarray*}\lefteqn{\bigg[F(U'(x)\tilde{Z}_{k}(u))\parallel\bar{\Theta}(W)(u) \parallel^2-  F(U'(x)\tilde{Z}(u))\parallel\bar{\Theta}(W)(u) \parallel^2\bigg]}\notag\\  &=&\bigg[F(U'(x)\tilde{Z}_{k}(u))\parallel\bar{\Theta}(W)(u) \parallel^2-F(U'(x)\tilde{Z}(u))\parallel\bar{\Theta}(W)(u) \parallel^2 \bigg]1_{\{\tau_{k}\leq u\}}\longrightarrow 0,\,\,\, \tilde{\mathbb{Q}}\,\,a.s.\end{eqnarray*} This is due to $\tau_{k}\uparrow T$\,\,\,$\mathbb{P}$\,\,a.s.\, (see Lemma \ref{L:1}) and $\tilde{\mathbb{Q}}\sim\mathbb{P}.$ 
In the light of (\ref{qwe1a}), (\ref{0op}) and (\ref{ris}), for fixed $u,$ the sequence $F(U'(x)\tilde{Z}_{k}(u))\parallel\bar{\Theta}(W)(u) \parallel^2$
indexed by $k$ is uniform integrable. Consequently 
\begin{equation}\label{P}
\mathbb{\tilde{E}}F(U'(x)\tilde{Z}_{k}(u))\parallel\bar{\Theta}(W)(u) \parallel^2\longrightarrow \mathbb{\tilde{E}}F(U'(x)\tilde{Z}(u))\parallel\bar{\Theta}(W)(u) \parallel^2.
\end{equation}  
Now \eqref{qwe1a}, \eqref{0op}, (\ref{ris}) and Lebesque's Dominated Convergence Theorem finish the proof of the Lemma.

\begin{flushright}
$\square$
\end{flushright} 

\begin{lemma}\label{L:4}
 $$\int_{0}^{T}\mathbb{\tilde{E}}(\lambda_{k_{n}}(t)\vert\mathcal{F}_{t})d\bar{W}(t)\longrightarrow\int_{0}^{T}\mathbb{\tilde{E}}(\lambda(t)\vert\mathcal{F}_{t})d\bar{W}(t),$$
for a subsequence $k_{n}\rightarrow\infty,$  $\mathbb{P}$ a.s.
\end{lemma}

\vspace{0.5cm}

\noindent {\sc  Proof:} 

\vspace{0.5cm}

 By It\^o's isometry it suffices to prove that
$$\mathbb{\tilde{E}}\int_{0}^{T}\|\mathbb{\tilde{E}}(\lambda_{k}(t)-\lambda(t)\vert\mathcal{F}_{t})\|^{2}\,dt \longrightarrow 0.$$
Moreover $$\mathbb{\tilde{E}}\int_{0}^{T}\|\mathbb{\tilde{E}}(\lambda_{k}(t)-\lambda(t)\vert\mathcal{F}_{t})\|^{2}\,dt\leq \int_{0}^{T}\mathbb{\tilde{E}}\|\lambda_{k}(t)-\lambda(t)\|^2\,dt,$$
due to Jensen's inequality. Hence we want to prove
\begin{equation}\label{023}
\int_{0}^{T}\mathbb{\tilde{E}}\|\lambda_{k}(t)-\lambda(t)\|^2\,dt\longrightarrow 0.
\end{equation}
Let us recall that
 $$\lambda_{k}(t,\omega)\triangleq\left[\int_{t}^{T}\mu(du,(\tilde{Z}_{k}(\omega),W(\omega)))\Phi_k(u,t,\omega)\right]g_{k}(t,(\tilde{Z}_{k}(\omega),W(\omega))),$$ and
 $$\lambda(t,\omega)\triangleq\left[\int_{t}^{T}\mu(du,(\tilde{Z}(\omega),W(\omega)))\Phi(u,t,\omega)\right]g(t,(\tilde{Z}(\omega),W(\omega))). $$
Because of $\tau_{k}\uparrow T$\,\,\,$\mathbb{P}$\,\,a.s.\, (see Lemma \ref{L:1}) and $\tilde{\mathbb{Q}}\sim\mathbb{P},$ for $t<T$ 
\begin{eqnarray*}\lefteqn{[g_{k}(t,(\tilde{Z}_{k}(\omega),W(\omega)))-g(t,(\tilde{Z}(\omega),W(\omega)))]}\\&&\mbox{}=[g_{k}(t,(\tilde{Z}_{k}(\omega),W(\omega)))-g(t,(\tilde{Z}(\omega),W(\omega)))]1_{\{\tau_{k}\leq t\}}\longrightarrow 0,\,\, \tilde{\mathbb{Q}}\,\,a.s. \end{eqnarray*} 
In order to prove the Lemma it suffices to show
\begin{equation}\label{show1}
\mathbb{\tilde{E}}\left\|\int_{t}^{T}\mu(du,(\tilde{Z}_{k},W))\Phi_k(u,t)-\int_{t}^{T}\mu(du,(\tilde{Z},W))\Phi(u,t)\right\|^{2}\rightarrow 0,
\end{equation} and
\begin{equation}\label{show2}
\mathbb{\tilde{E}}\,\|\lambda_{k}(t)\|^{2+\epsilon}\leq K_{1}
,\end{equation} for some $\epsilon>0$ and a constant $K_1$ independent of $k$ and $t.$ Indeed (\ref{show1}) give the almost sure convergence (up to a subsequence) of $\lambda_{k}$ to $\lambda.$ Moreover (\ref{show2}) implies the uniform convergence of $\|\lambda_{k}\|^{2}$, and also yields (\ref{023}) by Lebesque Dominated Convergence Theorem. 
 To proceed, we need some bounds on $\Phi_{k}(t,\cdot),$ and $\Phi(t,\cdot)$ independent of $k$ and $t.$  Cf (\ref{oop})
\begin{equation}\label{oop1}
 \sup_{s\leq v\leq t}\parallel\Phi^{2,j}(v,s)\parallel\leq e^{K_{2}(t-s)}.
\end{equation}
Moreover $\Phi_{k}^{2j}=\Phi^{2j}$ and $\Phi^{21}\equiv0.$ 
Furthermore we prove for $j=1,\cdots,n+1$
\begin{equation}\label{112}
\mathbb{\tilde{E}}\,[\,\sup_{0\leq t\leq T}|\Phi_{k}^{1j}(t,s)|]^{m}\leq K_3
,\end{equation}and

\begin{equation}\label{113}
\mathbb{\tilde{E}}\,[\,\sup_{0\leq t\leq T}|\Phi^{1j}(t,s)|]^{m}\leq K_4
,\end{equation} for some constants $K_3,$ $K_4,$ and $m>1.$  We prove (\ref{113}), and (\ref{112}) follows similarly, since  $|\phi'_k(x)|\leq1.$ 
One has ($\delta_{1j}$ is the Kronecker delta)
\begin{eqnarray*}
\Phi^{1,j}(t,s)&=&\int_{s}^{t}\left(\parallel\bar{\Theta}(W)(u)\parallel^2\Phi^{1,j}(u,s)+2\tilde{Z}(u)\bar{\Theta}^{T}(W)(u)[\bar{\Theta}'(W)\Phi^{2,j}(\cdot,s)](u)\right)\,du\\& & \mbox{} +\int_{s}^{t}\left(\bar{\Theta}^{T}(W)(u)\Phi^{1,j}(u,s)+\tilde{Z}(u)[\bar{\Theta}'(W)\Phi^{2,j}(\cdot,s)](u)\right)d\tilde{W}(u) +\delta_{1j}\\&=&\int_{s}^{t}\parallel\bar{\Theta}(W)(u)\parallel^2\Phi^{1,j}(u,s)\,du+\int_{s}^{t}\bar{\Theta}^{T}(W)(u)\Phi^{1,j}(u,s)\,d\tilde{W}(u)+A^{j}(t,s),
\end{eqnarray*} 
where
\begin{eqnarray}\notag
A^{j}(t,s) &\triangleq &\int_{s}^{t}2\tilde{Z}(u)\bar{\Theta}^{T}(W)(u)[\bar{\Theta}'(W)\Phi^{2,j}(\cdot,s)](u) \,du \\\notag
& & \mbox{}+\int_{s}^{t}\tilde{Z}(u)[\bar{\Theta}'(W)\Phi^{2,j}(\cdot,s)](u)\,d\tilde{W}(u) +\delta_{1j}.\label{09po}
\end{eqnarray}
Let $m>1,$ the inequality $$|a+b|^m\leq K_5(|a|^m+|b|^m),$$ implies
\begin{eqnarray}\label{ll}|\Phi^{1,j}(t,s)|^m &\leq & K_6\bigg[\left|\int_{s}^{t}\parallel\bar{\Theta}(W)(u)\parallel^2\Phi^{1,j}(u,s)\,du \right|^m\\\notag& & \mbox{}+\left|\int_{s}^{t}\bar{\Theta}^{T}(W)(u)\Phi^{1,j}(u,s)\,d\tilde{W}(u)\right|^m+|A^{j}(t,s)|^m\bigg], \notag
\end{eqnarray}
and
\begin{eqnarray*}
|A^{j}(t,s)|^m&\leq & K_7\bigg[\left|\int_{s}^{t}2\tilde{Z}(u)\bar{\Theta}^{T}(W)(u)[\bar{\Theta}'(W)\Phi^{2,j}(\cdot,s)](u)\,du\right|^m\\& & \mbox{}+\left|\int_{s}^{t}\tilde{Z}(u)[\bar{\Theta}'(W)\Phi^{2,j}(\cdot,s)](u)\,d\tilde{W}(u)\right|^m +\delta_{1j}\bigg]. 
\end{eqnarray*}
H\"older's inequality with  (\ref{as3}), (\ref{qwe1a}), (\ref{2}) and (\ref{oop1}) yield \begin{eqnarray*}
\lefteqn{\mathbb{\tilde{E}}\left|\int_{s}^{t}2\tilde{Z}(u)\bar{\Theta}^{T}(W)(u)[\bar{\Theta}'(W)\Phi^{2,j}(\cdot,s)](u) \,du\right|} \\& & \mbox{}\leq\mathbb{\tilde{E}}\left(\int_{s}^{t}\parallel\bar{\Theta}^{T}(W)(u) \parallel^{m}\,du \right)
\left[\mathbb{\tilde{E}}\int_{s}^{t} \parallel 2\tilde{Z}(u)[\bar{\Theta}'(W)\Phi^{2,j}(\cdot,s)](u) \parallel^{m'}\,du\right]^\frac{m}{m'}\leq K_{8},\end{eqnarray*}
with $m'$ the conjugate of $m,$ i.e., $\frac{1}{m}+\frac{1}{m'}=1.$ 

The Burkholder-Davis-Gundy inequality together with (\ref{as3}), (\ref{2}), (\ref{oop1})  and the above arguments show that

\begin{eqnarray*}\mathbb{\tilde{E}}\sup_{0\leq t\leq T}\left|\int_{s}^{t}\tilde{Z}(u)[\bar{\Theta}'(W)\Phi^{2,j}(\cdot,s)](u)\,d\tilde{W}(u)\right|^m &\leq&
\mathbb{\tilde{E}}\left|\int_{0}^{T}\parallel\tilde{Z}(u)[\bar{\Theta}'(W)\Phi^{2,j}(\cdot,s)](u) \parallel^{2}\,du\right|^{\frac{m}{2}}\\&\leq& K_{9}.
\end{eqnarray*}

Therefore
\begin{equation}\label{54}\mathbb{\tilde{E}}\,[\,\sup_{0\leq t\leq T}|A^{j}(t,s)|\,]^{m}\leq K_{10}.\end{equation}

H\"older's inequality and (\ref{qwe1a}) give
\begin{eqnarray}\label{nk}
\lefteqn{\mathbb{\tilde{E}}\left|\int_{s}^{t}\parallel\bar{\Theta}(W)(u)\parallel^2\Phi^{1,j}(u,s)\,du \right|^m}\\\notag & &\qquad\qquad\qquad\qquad\leq
\left[\int_{s}^{t}\mathbb{\tilde{E}}\,[\,\sup_{0\leq v\leq u}|\Phi^{1j}(v,s)|]^{m}\,du\right]\left[\mathbb{\tilde{E}}\int_{s}^{t}\parallel\bar{\Theta}(W)(u)\parallel^{2m'}\,du \right]^\frac{m}{m'}\\\notag& &\qquad\qquad\qquad\qquad\leq K_{11}\int_{0}^{T}\mathbb{\tilde{E}}\,[\,\sup_{0\leq v\leq u}|\Phi^{1j}(v,s)|]^{m}\,du.\notag
\end{eqnarray}
The Burkholder-Davis-Gundy inequality and the above arguments show
\begin{equation}\label{555}
\mathbb{\tilde{E}}\sup_{0\leq t\leq T}\left|\int_{s}^{t}\bar{\Theta}^{T}(W)(u)\Phi^{1,j}(u,s)d\tilde{W}(u)\right|^m\leq K_{12}\int_{0}^{T}\mathbb{\tilde{E}}\,[\,\sup_{0\leq v\leq u}|\Phi^{1j}(v,s)|]^{m}\,du.
\end{equation}

By combining (\ref{ll}), (\ref{54}), (\ref{nk}) and (\ref{555}) one gets
$$\mathbb{\tilde{E}}\,[\,\sup_{0\leq t\leq T}|\Phi^{1j}(t,s)|]^{m}\leq K_{13}+K_{14}\int_{0}^{T}\mathbb{\tilde{E}}\,[\,\sup_{0\leq v\leq u}|\Phi^{1j}(v,s)|]^{m}\,du .$$  Gronwall's inequality gives then 
\begin{equation}\label{zax}\mathbb{\tilde{E}}\,[\,\sup_{0\leq t\leq T}|\Phi^{1j}(t,s)|]^{m}\leq K_{15}. \end{equation}
On $[0,\tau_{k}]$ $\Phi_{k}=\Phi$ since it satisfies the same SDE. Thus for $j=1 ,\cdots, n+1$
\begin{eqnarray*}
\lefteqn{\left(\int_{t}^{T}\mu(du,(\tilde{Z}_{k},W))\Phi_k(u,t)-\int_{t}^{T}\mu(du,(\tilde{Z},W))\Phi(u,t)\right)_{j}}\\
& =&\int_{\tau_{k}}^{T}\left([F'(U'(x)\tilde{Z}_{k}(u))]\parallel\bar{\Theta}(W)(u)\parallel^{2}\Phi_{k}^{1,j}(u,t)\right)\,du\\&-&
\int_{\tau_{k}}^{T}\left([F'(U'(x)\tilde{Z}_{k}(u))]\parallel\bar{\Theta}(W)(u)\parallel^{2}\Phi^{1,j}(u,t)\right)\,du\\
&+&\int_{\tau_{k}}^{T}\left(2F(U'(x)\tilde{Z}_{k}(u))\bar{\Theta}^{T}(W)(u)[\bar{\Theta}'(W)\Phi_{k}^{2,j}(\cdot,t)](u)\right)\,du\\&-&
\int_{\tau_{k}}^{T}\left(2F(U'(x)\tilde{Z}(u)\bar{\Theta}^{T}(W)(u)[\bar{\Theta}'(W)\Phi^{2,j}(\cdot,t)](u)\right)\,du
\end{eqnarray*}
H\"older's inequality with (\ref{as3}), (\ref{qwe1a}), (\ref{0op}), (\ref{ris}), (\ref{2}), (\ref{oop1}), (\ref{zax}) and $\tau_{k}\uparrow T$\,\,\,$\tilde{Q}$\,\,a.s. prove (\ref{show1}).

If $|\cdot|$ denotes the matrix norm and $\varrho$ a nonnegative number, by  (\ref{qwe1a}), (\ref{ris}) and (\ref{2}),
$$\mathbb{\tilde{E}}\,|g_{k}(t,(\tilde{Z}_{k},W) |^{\varrho}<K,\quad\mbox{and}\quad\mathbb{\tilde{E}}\,|g(t,(\tilde{Z},W)|^{\varrho}<K,            $$
for a constant $K$ independent of $k$ and $t.$ This combined with H\"older's inequality and the above arguments prove (\ref{show2}).

\begin{flushright}
$\square$
\end{flushright} 
\vspace{0.5cm}

Lemmas \ref{L:3}, and \ref{L:4} in conjunction with the equation (\ref{09}) conclude the proof.

\begin{flushright}
$\square$
\end{flushright}

{\bf{Acknowledgements}}

\vspace{0.4cm}

The authors would like to thank Professor Steven E. Shreve for helpful discussions and comments.

\end{document}